\documentclass[a4paper]{amsart}
\usepackage{amssymb}
\usepackage{amsmath}
\usepackage{amsthm}
\usepackage{amsfonts}
\usepackage[T1]{fontenc}

\usepackage[colorlinks=true, linkcolor=blue, citecolor=red] {hyperref}
\usepackage{enumitem}
\usepackage{graphicx}
\usepackage{mathtools}
\allowdisplaybreaks

\numberwithin{equation}{section}

\newtheorem{theorem}{Theorem}
\newtheorem{lemma}[theorem]{Lemma}

\newtheorem{definition}[theorem]{Definition}
\newtheorem{remark}[theorem]{Remark}

\newcommand{\test}{\varphi}
\newcommand{\R}{\mathbb{R}}
\newcommand{\N}{\mathbb{N}}

\begin{document}

\title[Wellposedness for the spatially dependent OH Equation]{Wellposedness of the initial value problem for the Ostrovsky-Hunter Equation with spatially dependent flux}

\author[G. M. Coclite]{G. M. Coclite}
\address[Giuseppe Maria Coclite]
{\newline Dipartimento di Meccanica, Matematica e Management
\newline Politecnico di Bari
\newline Via E. Orabona 4
\newline 70125 Bari, Italy} 
\email[]{giuseppemaria.coclite@poliba.it}

\author[N. Chatterjee]{N. Chatterjee}
\address[Neelabja Chatterjee]
{\newline Department of mathematics
\newline University of Oslo
\newline P.O. Box 1053,  Blindern
\newline N--0316 Oslo, Norway} 
\email[]{neelabjc@math.uio.no}

\author[N. H. Risebro]{N. H. Risebro}
\address[Nils Henrik Risebro]
{\newline Department of mathematics
\newline University of Oslo
\newline P.O. Box 1053,  Blindern
\newline N--0316 Oslo, Norway} 
\email[]{nilsr@math.uio.no}

\date{\today}


\thanks{This project has received funding from the European Union's
Horizon 2020 research and innovation programme under the Marie Sk\l{}odowska-Curie grant
agreement No 642768 for NC's work and GMC's work was supported by the Research Council of Norway (project 250674/F20). GMC is  member of the Gruppo Nazionale per l'Analisi Matematica, la Probabilita' e le loro Applicazioni (GNAMPA) of the Istituto Nazionale di Alta Matematica (INdAM). This work was initiated while GMC visited Department of Mathematics  at the University of Oslo. He is grateful for Department's financial support and excellent working conditions.}

\begin{abstract}
  In this paper we study the Ostrovsky-Hunter equation for the case 
  where the flux function $f(x,u)$ may depend on the spatial variable 
  with certain smoothness. Our main results are that if the flux 
  function is smooth enough (namely $f_{x}(x,u)$ is uniformly Lipschitz 
  locally in $u$ and $f_{u}(x,u)$ is uniformly bounded), then there exists
  a unique entropy solution. To show the existence, after proving some 
  \textit{a priori} estimates we have used the method of 
  \textit{compensated compactness} and to prove the uniqueness we have 
  employed the method of \textit{doubling of variables}.
\end{abstract}

\maketitle

\section{Introduction}

To model numerous physical phenomena such as the propagation of undular bores in shallow water, 
the flow of liquids containing gas bubbles, the propagation of waves in an elastic tube filled with
a viscous fluid, weakly nonlinear plasma waves with certain dissipative effects \textit{etc.} the 
following nonlinear evolution equation, known as \textit{Korteweg-deVries-Burgers equation}
\begin{equation}\label{eq:KdvBurger}
u_{t} + \Big(f(u)\Big)_{x} - \alpha u_{xx} - \beta u_{xxx} = 0, 
\quad \alpha, \beta \in \R, \quad f(u) = \frac{u^{2}}{2},
\end{equation}
has been extensively investigated in the recent years (see  \cite{Gallego, KdV, Shu} and references therein). Also considering the effects 
of background rotation through the Coriolis force ($\kappa$ being the force parameter and $C_{0}$ is 
the linear shallow water speed) \eqref{eq:KdvBurger} takes the following form
\begin{equation}\label{eq:KdVBurgerI}
\Big(u_{t} + \Big(f(u)\Big)_{x} - \alpha u_{xx} - \beta u_{xxx} \Big)_{x} = \gamma u, \quad \gamma = \frac{\kappa^{2}}{2C_{0}}> 0.
\end{equation}
To model small-amplitude long waves in a rotating fluid of finite depth \cite{Ostrovsky} and to study long internal waves in a rotating fluid \cite{HunterI} both the viscous dissipation term and the high-frequency dispersion term has to be dropped, \textit{i.e.} $\alpha, \beta = 0$; which leads to
\begin{equation}\label{eq:OHI}
\Big(u_{t} + \Big(f(u)\Big)_{x} \Big)_{x} = \gamma u,
\end{equation}
which is known as the \textit{Ostrovsky-Hunter equation}, as Ostrovsky also independently derived 
them \cite{Ostrovsky}. This equation is also used to model high frequency waves in a relaxing medium 
\cite{VakhnenkoI}. In the cases described above the flux is considered to be of Burger's type, 
\textit{i.e.} $f(u) = \frac{u^{2}}{2}$.

Also by including the effects of background rotation in the shallow water equation, and then using 
singular perturbation methods \eqref{eq:OHI} has been derived previously (see \cite{RuvoThesis}, \cite{HunterII}). In the recent years enormous amount of research has been carried out investigating \eqref{eq:OHI}. Among those works in \cite{Ostrovsky}, \cite{ParkesI}, \cite{Stepanyants} the equation \eqref{eq:OHI} is also known as the reduced Ostrovsky equation, in \cite{HunterI} it is called short wave equation, whereas in \cite{Boutet}, \cite{Brunelli} \eqref{eq:OHI} is known as Ostrovsky-Vakhnenko equation and as Vakhnenko equation in \cite{VakhnenkoIII}. Moreover, the equation \eqref{eq:KdVBurgerI} is used to model ultra short light pulses in silica optical fibres (see \cite{AVB}, \cite{Liu}), in which case $f(u) = - \frac{1}{6} u^{3}$. In this case equation \eqref{eq:KdVBurgerI} is sometimes referred to as the \textit{short-pulse-equation}. 

In his seminal paper \cite{HunterI}, Hunter showed the connection between the KdV equation 
\eqref{eq:KdVBurgerI} and the Ostrovsky-Hunter equation \eqref{eq:OHI} as the no-rotation and no-long wave dispersion limits of the same equation. When the oceanic waves approach shore, the waves usually propagate through a background with varying properties. It is natural to expect the linear phase speed of the wave which encoded in the flux function, in such a variable medium, should have a \textit{spatial dependecy}. In the context of KdV equation, Johnson \cite{Johnson} for water waves and Grimshaw \cite{GrimshawI} for internal waves derived the variable coefficient equation (see also 
\cite{GrimshawII} for a detailed review). Motivated by this, it is immediate to pose the question of design and analysis of numerical scheme for the Ostrovsky-Hunter equation with a spatial dependency in the flux function. In \cite{Neelabja}, we investigated the \textit{spatially dependent} Ostrovsky-Hunter equation in the fully-discretized setting to prove convergence of the corresponding numerical method to the unique entropy solution and we proved its \textit{order of convergence}. Whereas in this paper, our aim is to establish well-posedness of the problem in continuous set up.

The results obtained in this paper are the following. If the function $f_{x}(x,u)$ is uniformly Lipschitz continuous locally in $u$, the function $f_{u}(x,u)$ to be uniformly locally bounded, and the initial data are square integrable and satisfy zero-mean condition, then there exists an entropy solution via method of compensated compactness. Furthermore, for two entropy solutions $u$ and $v$, with initial data $u_{0}$ and $v_{0}$ respectively, we establish the following estimate 
\begin{equation*}
  \|u(\cdot, t) - v(\cdot, t)\|_{L^{1}((0,R))} \leq e^{Ct}\|u_{0}(\cdot) - v_{0}(\cdot)\|_{L^{1}((0,R+Lt))}, \quad
\end{equation*}
for some constants $C, L\;\mathrm{and}\;R$ to be specified later. The rest of this paper is organized as follows. In Section \ref{sec:Prelim} we give detailed descriptions of the notations used, the precise assumptions of the regularity of the flux funxtion and the initial data. Also, apart from stating our main result as a theorem, we state the definition of entropy solution to be used. In Section \ref{sec:Estimates} we prove few useful \textit{a priori} estimates (namely energy estimate and $L^{\infty}_{\mathrm{loc}}$ bound) for the purpose of compensated compactness technique. In Section \ref{sec:MainResult} we first state the two results due to Murat and Tartar in the form of two lemmas, using which we will employ a compensated compactness argument to show the existence of entropy solution of the equation under consideration. Moreover, we establish an $L^{1}$ contraction type estimate mentioned above using the technique of doubling of variables.

\section{Preliminaries and Notation}\label{sec:Prelim}

Throughout this paper $u(x,t)$ is the conserved quantity and $f$ is the flux which is
dependent on the spatial variable $x$ and $u(x,t)$, denoted by $f(x,u(x,t))$. For notational consistency, we mention the following chain rule keeping the notation of $f_{x}(x,u) \neq f(x,u)_{x}$
\begin{equation*}
\begin{cases}
 \frac{\partial f}{\partial u}(x, u) = f_{u}(x, u), \; \\
 \partial_{x}f(x,u) = f(x,u)_{x} = f_{u}(x, u)\frac{\partial u}{\partial x} + f_{x}(x, u), \; \\
 \frac{\partial u}{\partial t} (x,t) = u_{t}(x,t).
 \end{cases}
\end{equation*}
We are interested in the initial boundary value problem for \eqref{eq:OHI}, but with \textit{spatially dependent} flux, and hence we augment the equation with the initial datum
\begin{equation}\label{eq:InitialDatum}
u(x,0) = u_{0}(x), \; \mathrm{for} \; x > 0.
\end{equation}
Keeping that in mind, following the works of \cite{GiuseppeV}, \cite{GiuseppeI} and \cite{LiuI} integrating the equation \eqref{eq:OHI} on the interval $(0, x)$ we get the integro-differential formulation of the problem under consideration and setting $\gamma =1$,
\begin{equation}\label{eq:OHII}
\begin{cases}
u_{t} + f(x,u)_{x} = \int\displaylimits_{0}^{x} u(y,t) dy, \; t>0, \; x > 0, \\
u(x,0) = u_{0}(x), \qquad \qquad \qquad \qquad \; \;  x > 0, \\
u(0,t) = 0, \qquad \qquad \qquad \qquad \qquad \; \; \, t > 0.\\
\end{cases}
\end{equation}
Denoting $P[t,x] := \int\displaylimits_{0}^{x} u(y,t) dy$, we get the following equivalent formulation
\begin{equation}\label{eq:OHIII}
\begin{cases}
u_{t} + f(x,u)_{x} = P[t,x], \; t>0, \; x > 0, \\
P_{x} = u(x,t), \; \; \qquad \; \qquad t>0, \; x > 0, \\
P(t, 0) = u(0,t) = 0, \; \; \quad \qquad \; \; \, t>0, \\
u(x,0) = u_{0}(x), \qquad \qquad \qquad \; \; x > 0.
\end{cases}
\end{equation}
For the initial datum, we assume the following \textit{zero-mean condition} and regularity assumption respectively
\begin{equation}\label{eq:ZeroMeanCondition}
\int\displaylimits_{0}^{\infty} u_{0}(x) dx = 0,
\end{equation}
and
\begin{equation}\label{eq:RegularityAssumption}
u_{0}(x) \in L^{2}(\R_{+}) \cap L^{\infty}_{\mathrm{loc}}(\R_{+}),
\end{equation}
where $\R_{+} := (0,\infty)$ will be denoting the unbounded positive half line throughout the paper. Similarly $\Pi$ will be used to denote $(0,\infty)^{2}$. Also the flux $f$ is assumed to satisfy the following hypothesis:
\begin{enumerate}[label=(A\arabic*)]
 \item \label{cond:genuinenonlinearity} $f(x,\cdot)$ is \textit{genuinely nonlinear}, i.e. $f_{uu}(x,u) \neq 0$ for a.e. $(x,u) \in 
 \R_{+} \times \R$ and $\displaystyle{\lim_{x \to \infty}}\partial_{x}f(x,u) = \displaystyle{\lim_{x \to \infty}} f(x,u) = \displaystyle{\lim_{x \to 0}} f(x,u) = \displaystyle{\lim_{x \to 0}}\partial_{x}f(x,u) = 0$ for all $u$,
 \item \label{cond:boundedderivative} $\exists$ a constant $C > 0$ such that $|f_{xu}(x,u)| \leq C$ and $|f_{x}(x,u)| \leq C|u|$ for all $u$,
 \item \label{cond:xLipschitz} $\exists$ a constant $L_{1} > 0$ such that $|f_{x}(x,u) - f_{x}(x,v)| \leq L_{1}|u - v|$, for all $u, v$,
 \item \label{cond:uLipschitz} $\exists$ a constant $L > 0$ such that $|f_{u}(x,u)| \leq L$, for all $u$.
\end{enumerate}
Even if the initial data is smooth enough, solutions of \eqref{eq:OHIII} generically develop discontinuities. Hence solutions must be considered in the weak sense. A function $u$ is a \emph{weak solution} of \eqref{eq:OHIII} if
\begin{equation} \label{eq:WeakSoln}
  \begin{aligned}
    \iint\displaylimits_{\Pi} u
    \test_{t} & + f(x,u) \test_{x} + P[t,x] \test \, dx dt + \int\displaylimits_{\R_{+}} u_{0}(x) \test(x,0) \,dx  = 0,
  \end{aligned}
\end{equation}
for all test functions $\test=\test(x,t) \in C_{c}^{\infty}(\Pi)$. Moreover, from \eqref{eq:OHIII} we have that
\begin{equation}\label{eq:regularization}
 u \in L^{\infty}_{\mathrm{loc}}(\Pi) \Rightarrow P \in L^{\infty}_{\mathrm{loc}}(\R_{+}; W^{1, \infty}_{\mathrm{loc}}\R_{+}).
\end{equation}
Following \cite{Bardos} we define entropy solutions as
\begin{definition}[Entropy Solution]\label{def:EntropySoln}
 We say that $u \in L^{\infty}_{\mathrm{loc}}\Big(\Pi\Big)$ is an \textit{entropy solution} of the initial boundary value problem \eqref{eq:OHIII}, if 
 \begin{itemize}
   \item $u$ satisfies \eqref{eq:WeakSoln} ; 
   
   \item for every smooth, non negative test function $\phi \in C^{2}_{c}(\Pi)$ and $c \in \R$
   \begin{align}\label{eq:StrongTraceEntropy}
   & \iint\displaylimits_{\Pi} \Big(|u-c| \partial_{t} \phi + \mathrm{sign}(u-c)\Big(f(x,u) - f(x,c)\Big) \partial_{x} \phi - \mathrm{sign}(u-c)f_{x}(x,c) \phi \Big) dt dx \nonumber \\ & \qquad + \iint\displaylimits_{\Pi} \mathrm{sign}(u-c) P \phi dt dx - \int\displaylimits_{\R_{+}} \mathrm{sign}(c)\Big(f(0,u^{\tau}_{0}) - f(0,c)\Big) \phi(0,t) dt \nonumber \\ & \qquad + \int\displaylimits_{\R_{+}} |u_{0}(x) - c| \phi(x,0) dx \geq 0 .
   \end{align}
 \end{itemize}
\end{definition}
As an immediate consequence of \eqref{eq:regularization} if a map $u \in L^{\infty}_{\mathrm{loc}}(\Pi)$ satisfies the following equivalent entropy inequality for every convex entropy/entropy flux pair (i.e. for $\eta \in C^{2}(\R), q(x,u): = \int_{0}^{u}\eta'(v)f_{v}(x,v)dv$)
\begin{equation}\label{eq:entropyinequality}
\partial_{t} \eta(u) + \partial_{x} q(x,u) + \eta'(u)f_{x}(x,u) - q_{x}(x,u) - \eta'(u)P \leq 0, 
\end{equation}
in the sense of distributions, then by Theorem $1.1$ of \cite{GiuseppeY} on the boundary $x=0$ strong trace $u^{\tau}_{0}$ exists. By a standard approximation argument equivalently any \textit{convex entropy/entropy flux pair} $(\eta,q)$ in \eqref{eq:entropyinequality} can be replaced by Kru\v{z}kov entropy pair namely for $c \in \R$, $\eta(u) = |u-c|$ and $q(x,u) = \int\displaylimits_{0}^{u} \, \mathrm{sign}(u-c) \, f_{v}(x,v) \, dv$. 
The main result of this paper is the following theorem.
\begin{theorem}\label{thm:wellposedness}
Assuming \eqref{eq:ZeroMeanCondition} and \eqref{eq:RegularityAssumption}, the Cauchy problem \eqref{eq:OHII}, or equivalently \eqref{eq:OHIII} possesses a unique entropy solution $u$ in the sense of Definition \ref{def:EntropySoln}. Moreover, if $u$ and $v$ are two entropy solutions of \eqref{eq:OHII}, or equivalently \eqref{eq:OHIII} in the sense of Definition \ref{def:EntropySoln}, the following estimate holds for a given $0 < t \leq T$ 
 \begin{equation}\label{eq:contractivity}
  \|u(\cdot, t) - v(\cdot, t)\|_{L^{1}((0,R))} \leq e^{Ct}\|u_{0}(\cdot) - v_{0}(\cdot)\|_{L^{1}((0,R+Lt))} \quad
 \end{equation}
for almost every $T > t > 0$, $R > 0$ and $L > 0$ being the bound $|f_{u}(x,u)| \leq L$, where the constant $C$ depends on $T$, $R$,
and $L$.
\end{theorem}

Before proceeding to prove this theorem, it is worth mentioning that Coclite \textit{et al.} \cite{GiuseppeV}, \cite{GiuseppeI} have showed the well-posedness of the initial-boundary value problem and the Cauchy problem for the Ostrovsky-Hunter Equation \eqref{eq:OHII}, but
without any spatial dependency in the flux. Throughout the next section we will extend their results following the papers cited just above. 

\section{A-Priori Estimates}\label{sec:Estimates}
The existence argument is based on passing to the limit in the following vanishing viscosity approximation of \eqref{eq:OHIII} (see \cite{GiuseppeIV}). Fix a small number $\epsilon > 0$, and let $u_{\epsilon} = u_{\epsilon}(x,t)$ be the unique classical solution of the following problem 
\begin{equation}\label{eq:OHIV}
 \begin{cases}
 \partial_{t} u_{\epsilon} + \partial_{x} f(x,u_{\epsilon}) = P_{\epsilon} + \epsilon \partial_{xx}^{2}u_{\epsilon}, \qquad &t>0, x > 0 \\
 \partial_{x}P_{\epsilon} = u_{\epsilon},   \qquad &t>0, x > 0 \\
 P_{\epsilon}(t,0) = u_{\epsilon}(0,t) =0, \qquad &t>0 \\
 u_{\epsilon}(x,0) = u_{\epsilon,0}(x), \qquad &x > 0,
 \end{cases}
\end{equation}
where $u_{\epsilon,0}$ is a $C^{\infty}(\R_{+})$ approximation of $u_{0}$ such that
\begin{equation}\label{eq:ViscousRegularity}
 \|u_{\epsilon,0}\|_{L^{2}(\R_{+})} \leq \|u_{0}\|_{L^{2}(\R_{+})}, \quad \int\displaylimits_{\R_{+}} u_{\epsilon,0}(x) dx = 0,
\end{equation}
and on the viscous source term for $x > 0, \; P_{\epsilon,0}(x) := \int\displaylimits_{0}^{x} u_{\epsilon,0}(y) dy$ we assume that
\begin{equation}\label{eq:ViscousSourceRegularity}
\begin{cases}
\|P_{\epsilon,0}\|^{2}_{L^{2}(\R_{+})} = \int\displaylimits_{\R_{+}} \Big(\int\displaylimits_{0}^{x} u_{\epsilon,0}(y) dy\Big)^{2} dx < \infty, \\
\int\displaylimits_{\R_{+}} P_{\epsilon,0}(x) dx = \int\displaylimits_{\R_{+}} \Big(\int\displaylimits_{0}^{x} u_{\epsilon,0}(y) dy\Big) dx = 0.
\end{cases}
\end{equation}
Clearly, \eqref{eq:OHIV} is equivalent to the integro-differential problem
\begin{equation}\label{eq:OHV}
 \begin{cases}
    \partial_{t} u_{\epsilon} + \partial_{x} f(x,u_{\epsilon}) = \int\displaylimits_{0}^{x}u_{\epsilon}(y,t)dy  + \epsilon \partial_{xx}^{2}u_{\epsilon}
    \qquad &t>0, x > 0, \\    
    u_{\epsilon}(x,0) = u_{\epsilon,0}(x) \qquad &x > 0.
 \end{cases}
\end{equation}
The existence of such solutions can be obtained by fixing a small number $1 > \delta > 0$ and considering the further approximation of \eqref{eq:OHV} (see for the whole real line \cite{GiuseppeI, GiuseppeV}; for the half line \cite{GiuseppeII, GiuseppeIII, GiuseppeIV} and the references therein). We are going to use the following result from previous works of Coclite \textit{et al.} (see \cite{GiuseppeI, GiuseppeX} and references therein).
\begin{theorem}\label{thm:AvailableResult}
Let $T>0$. Assuming that conditions \eqref{eq:ViscousRegularity} and \eqref{eq:ViscousSourceRegularity} hold, there exists a unique classical solution to the Cauchy problem of \eqref{eq:OHV} such that
\begin{equation}
\begin{cases}
u_{\epsilon} \in L^{\infty}_{\mathrm{loc}}\Big((0,T) \times \R_{+} \Big) \cap C \Big((0,T); H^{l}(\R_{+})\Big), \; \text{for all} \; l \in \N, \\
P_{\epsilon} \in L^{\infty}_{\mathrm{loc}}\Big((0,T) \times \R_{+} \Big) \cap L^{2}\Big((0,T) \times \R_{+} \Big), \\
\int\displaylimits_{0}^{\infty} u_{\epsilon}(x,t) dx = 0, \; t \geq 0.
\end{cases}
\end{equation}
\end{theorem}
Now let us prove some \textit{a priori} estimates on $u_{\epsilon}$.
\begin{lemma}\label{lem:aprioriequiv}
We have the equivalence of following two equalities
\begin{align} 
\begin{split}\label{eq:ViscousZeroMeanCondition}
 & \int\displaylimits_{\R_{+}}  u_{\epsilon}(x,t) dx = 0; \; t \geq 0,
 \end{split}\\
 & \|u_{\epsilon}(\cdot,t)\|_{L^{2}(\R_{+})}^{2} + 2 \epsilon \int\displaylimits_{0}^{t}\|\partial_{x}u_{\epsilon}(\cdot,s)\|_{L^{2}(\R_{+})}^{2}ds \nonumber \\ & = 
  \|u_{\epsilon,0}(\cdot)\|_{L^{2}(\R_{+})}^{2} + 2\int\displaylimits_{0}^{t}\Big[\int\displaylimits_{\R_{+}}[\int\displaylimits_{0}^{u_{\epsilon}}vf_{xv}(x,v)dv - u_{\epsilon} f_{x}(x,u_{\epsilon})]dx \Big] ds; \; t>0. \label{eq:EnergyEquality}
 \end{align}
\begin{proof}
Let $t>0$.  First we will prove that \eqref{eq:ViscousZeroMeanCondition} implies \eqref{eq:EnergyEquality}. Multiplying equation \eqref{eq:OHV} by $u_{\epsilon}(x,t)$ we get
\begin{align}\label{eq:aprioricalculationI}
 u_{\epsilon} \partial_{t}u_{\epsilon} + u_{\epsilon} f_{u}(x,u_{\epsilon}) (\partial_{x} u_{\epsilon}) + u_{\epsilon} f_{x}(x, u_{\epsilon}) &= 
  u_{\epsilon}\int\displaylimits_{0}^{x}u_{\epsilon}(y,t)dy + \epsilon u_{\epsilon} \partial_{xx}^{2} u_{\epsilon} .
\end{align}
In this equality, consider the term $ u_{\epsilon}\int\displaylimits_{0}^{x}u_{\epsilon}(y,t)dy$. We are going to show that after integration this term vanishes. Clearly the equation \eqref{eq:ViscousZeroMeanCondition} implies that
\begin{align}\label{eq:aprioricalculationII}
 \int\displaylimits_{\R_{+}} u_{\epsilon}(x,t) \Big[\int\displaylimits_{0}^{x} u_{\epsilon}(y,t) dy\Big] dx &=  \int\displaylimits_{\R_{+}} P_{\epsilon}( \partial_{x} P_{\epsilon}) dx =  \int\displaylimits_{\R_{+}} \partial_{x}[\frac{1}{2} P_{\epsilon}^{2}] dx =0.
\end{align}
Using $\eta(u) := \frac{1}{2}u^{2}$ into \eqref{eq:aprioricalculationI} we get
\begin{equation*}
 [u_{\epsilon}^{2}(x,t)]_{t} + 2 [q(x, u_{\epsilon})]_{x} - 2 \epsilon u_{\epsilon} \partial_{xx}^{2} u_{\epsilon} = 
 2 \Big[ \int_{0}^{u_{\epsilon}} v f_{xv}(x,v) dv -u_{\epsilon} f_{x}(x,u_{\epsilon}) \Big] + \partial_{x} [\frac{1}{2} P_{\epsilon}^{2}].
\end{equation*}
Integrating this expression over $\R_{+}$ and invoking \eqref{eq:aprioricalculationII} we get 
\begin{equation}\label{eq:aprioricalculationIII}
 \frac{d}{dt}\|u_{\epsilon}(t,\cdot)\|^{2}_{L^{2}(\R_{+})} + 2 \epsilon \|\partial_{x}u_{\epsilon}\|^{2}_{L^{2}(\R_{+})} = 2 \int\displaylimits_{\R_{+}} 
 \Big[ \int\displaylimits_{0}^{u_{\epsilon}} v f_{xv}(x,v) dv -u_{\epsilon} f_{x}(x,u_{\epsilon})  \Big] dx.
\end{equation}
And finally integrating over the $(0,t)$ we obtain
\begin{align*}
\begin{split}
 &\|u_{\epsilon}(\cdot,t)\|^{2}_{L^{2}(\R_{+})} + 2 \epsilon \int\displaylimits_{0}^{t}\|\partial_{x}u_{\epsilon}(\cdot,s)\|^{2}_{L^{2}(\R_{+})} ds \\
 & \qquad \qquad = \|u_{\epsilon,0}(\cdot)\|_{L^{2}(\R_{+})}^{2} + 2\int\displaylimits_{0}^{t}\Big[\int\displaylimits_{\R_{+}}[\int\displaylimits_{0}^{u_{\epsilon}}vf_{xv}(x,v)dv - u_{\epsilon} f_{x}(x,u_{\epsilon})]dx \Big] ds, 
 \end{split}
\end{align*}
which proves \eqref{eq:EnergyEquality}. Now we are going to prove the opposite implication. Assume that $\int\displaylimits_{\R_{+}} u_{\epsilon}(x,t) dx \neq 0$ for some $t > 0$, which implies
\begin{equation*}
P_{\epsilon}^{2}(t,+\infty) = \Big(\int\displaylimits_{\R_{+}} u_{\epsilon}(x,t) dx\Big)^{2} \neq 0,
\end{equation*}
which results in 
\begin{equation*}
 \frac{d}{dt}\|u_{\epsilon}(t,\cdot)\|^{2}_{L^{2}(\R_{+})} + 2 \epsilon \|\partial_{x}u_{\epsilon}\|^{2}_{L^{2}(\R_{+})} \neq 2 \int\displaylimits_{\R_{+}} 
 \Big[ \int\displaylimits_{0}^{u_{\epsilon}} v f_{xv}(x,v) dv -u_{\epsilon} f_{x}(x,u_{\epsilon})  \Big] dx,
\end{equation*}
ultimately contradicting our assumption \eqref{eq:EnergyEquality}.
This concludes the proof.
\end{proof}
\end{lemma}

\begin{lemma}\label{lm:zeromeanandenergybound}
For each $t \geq 0$, \eqref{eq:ViscousZeroMeanCondition} holds. In particular we have that for a constant $C>0$ coming from \ref{cond:boundedderivative}, independent of $\epsilon$
 \begin{align}\label{eq:EnergyEstimate}
  \|u_{\epsilon}(\cdot,t)\|^{2}_{L^{2}(\R_{+})} + 2 \epsilon \int\displaylimits_{0}^{t}\|\partial_{x}u_{\epsilon}(\cdot,s)\|^{2}_{L^{2}(\R_{+})} ds 
  & \leq \|u_{0}\|^{2}_{L^{2}(\R_{+})} \nonumber \\ & + \hat{C} \int\displaylimits_{0}^{t} \|u_{\epsilon}(\cdot,s)\|^{2}_{L^{2}(\R_{+})} ds
 \end{align}
where $\hat{C}$ is any constant greater than $C$.
 
\begin{proof}
  From the equation \eqref{eq:OHV} we have
  \begin{equation*}
   \partial_{x}\big(\partial_{t}u_{\epsilon} + \partial_{x}f(x,u_{\epsilon}) - \epsilon \partial_{xx}^{2} u_{\epsilon}\big) = u_{\epsilon}.
  \end{equation*}
Integrating both sides with respect to $x$ we get
  \begin{equation*}
  \partial_{t}u_{\epsilon} + \partial_{x}f(x,u_{\epsilon}) - \epsilon \partial_{xx}^{2} u_{\epsilon}\Big|_{0}^{\infty}  = \int\displaylimits_{\R_{+}} u_{\epsilon} dx.
  \end{equation*}
Observe that from $u_{\epsilon}(0,t) = 0$ of \eqref{eq:OHIV} we have $\partial_{t}u_{\epsilon}(0,t) = 0$ which, due to \ref{cond:genuinenonlinearity}
\begin{equation}\label{eq:initialdata}
\epsilon \partial^{2}_{xx} u_{\epsilon}(0,t) = \partial_{t}u_{\epsilon}(0,t) + \partial_{x}f(x,u_{\epsilon}) \Big|_{x = 0} - \int\displaylimits_{0}^{0} u_{\epsilon}(y,t) dy = 0.
\end{equation}
Invoking the property \ref{cond:genuinenonlinearity}, \eqref{eq:initialdata} and the smoothness of $u_{\epsilon}(x,t)$ from Theorem \ref{thm:AvailableResult} we can conclude $\int\displaylimits_{\R} u_{\epsilon}(x,t) dx = 0$, which proves \eqref{eq:ViscousZeroMeanCondition}. So by Lemma \ref{lm:zeromeanandenergybound} the relation \eqref{eq:EnergyEquality} holds.
To estimate the last term of the relation \eqref{eq:EnergyEquality} due to our assumption \ref{cond:boundedderivative} for any constant $\hat{C} \geq C$ we get
 \begin{equation}\label{eq:aprioricalculationIV}
   \Big|\int\displaylimits_{0}^{t}\Big[\int\displaylimits_{\R_{+}}[\int\displaylimits_{0}^{u_{\epsilon}}vf_{xv}(x,v)dv - f_{x}u_{\epsilon}(s,x)]dx \Big] ds \Big| \leq \hat{C} \int\displaylimits_{0}^{t} \|u_{\epsilon}(s, \cdot)\|^{2}_{L^{2}(\R_{+})} ds.
  \end{equation}
Consequently in \eqref{eq:EnergyEquality} inserting \eqref{eq:ViscousRegularity} and \eqref{eq:aprioricalculationIV} we have:
  \begin{align*}
    \|u_{\epsilon}(\cdot, t)\|_{L^{2}(\R_{+})}^{2} & + 2 \epsilon \int\displaylimits_{0}^{t}\|\partial_{x}u_{\epsilon}(\cdot, s)\|_{L^{2}(\R_{+})}^{2}ds \nonumber \\ \qquad &\leq \|u_{\epsilon,0}\|_{L^{2}(\R_{+})}^{2} + \hat{C} \int\displaylimits_{0}^{t} \|u_{\epsilon}(\cdot, s)\|^{2}_{L^{2}(\R_{+})} ds \\ \qquad & \leq \|u_{0}\|_{L^{2}(\R_{+})}^{2} + \hat{C} \int\displaylimits_{0}^{t} \|u_{\epsilon}(\cdot, s)\|^{2}_{L^{2}(\R_{+})} ds,
  \end{align*}
which concludes the proof.
 \end{proof}
\end{lemma}
\begin{remark}
It follows from \eqref{eq:EnergyEstimate} that
\begin{align*}
   \|u_{\epsilon}(\cdot, t)\|_{L^{2}(\R_{+})}^{2} &\leq \|u_{\epsilon}(\cdot, t)\|_{L^{2}(\R_{+})}^{2} 
   + 2 \epsilon \int\displaylimits_{0}^{t}\|\partial_{x}u_{\epsilon}(\cdot, s)\|_{L^{2}(\R_{+})}^{2}ds \\ &\leq
  \|u_{0}\|_{L^{2}(\R_{+})}^{2} + \hat{C} \int\displaylimits_{0}^{t} \|u_{\epsilon}(\cdot, s)\|^{2}_{L^{2}(\R_{+})} ds.
  \end{align*} 
 Thus by an application of Gronwall's inequality, we have
 \begin{equation}\label{eq:energybound}
  \|u_{\epsilon}(\cdot, t)\|_{L^{2}(\R_{+})} \leq e^{\hat{C}t} \|u_{0}\|_{L^{2}(\R_{+})}.
 \end{equation}
\end{remark}
\begin{lemma}\label{lm:LocallyBoundLemma}
The family
  \begin{equation}\label{eq:uLinfinitylocbound}
   \{u_{\epsilon}\}_{\epsilon > 0} \quad \text{is bounded in} \quad L^{\infty}_{\mathrm{loc}}(\Pi). 
  \end{equation}
  And consequently the family
  \begin{equation}\label{eq:PLinfinitylocbound}
   \{P_{\epsilon}\}_{\epsilon > 0} \quad \text{is bounded in} \quad L^{\infty}_{\mathrm{loc}}(\Pi).
  \end{equation}
\begin{proof}
By H\"{o}lder inequality we have the following estimate
  \begin{align*}
   \partial_{t}u_{\epsilon} + \partial_{x} f(x,u_{\epsilon}) - \epsilon \partial_{xx}^{2} u_{\epsilon} &= 
  \int\displaylimits_{0}^{x}u_{\epsilon}(t,y) dy \leq \Big|\int\displaylimits_{0}^{x}u_{\epsilon}(t,y) dy\Big| \\
   &\leq \int\displaylimits_{0}^{x}|u_{\epsilon}(t,y)| dy, \; \; \text{by H\"{o}lder's inequality,} \\
   &\leq \sqrt{x} \|u_{\epsilon}(t, \cdot)\|_{L^{2}(\R_{+})}, \; \; \text{using} \; \eqref{eq:energybound}) \\ 
   &\leq \sqrt{x} e^{\hat{C}t} ||u_{0}||_{L^{2}(\R_{+})}.
  \end{align*}
Now assume $v_{\epsilon}$ and $w_{\epsilon}$ be the solutions of the following equations respectively  
  \begin{equation}
   \begin{cases}
      \partial_{t}v_{\epsilon} + \partial_{x}f(x,v_{\epsilon}) =\|u_{0}\|_{L^{2}(\R_{+})} \sqrt{x} + \epsilon \partial_{xx}^{2}v_{\epsilon}, \quad t > 0, x > 0, \\
      v_{\epsilon}(0,x) = u_{\epsilon, 0}(x), \qquad \qquad \qquad \quad \qquad \qquad \qquad \qquad x > 0,
   \end{cases}
  \end{equation}  
  \begin{equation}
   \begin{cases}
      \partial_{t}w_{\epsilon} + \partial_{x}f(x,w_{\epsilon}) = - ||u_{0}||_{L^{2}(\R_{+})} \sqrt{x} + \epsilon \partial_{xx}^{2}w_{\epsilon}, \quad t > 0, x > 0, \\
      w_{\epsilon}(0,x) = u_{\epsilon, 0}(x), \qquad \qquad \qquad \qquad \qquad \qquad \qquad \qquad  \; x > 0.
   \end{cases}
  \end{equation}
  Then $u_{\epsilon}, v_{\epsilon},$ and $w_{\epsilon}$ are respectively a solution, a supersolution, and a subsolution of the parabolic problem \eqref{eq:OHV}. Following \cite[Theorem $9$, Chapter $2$]{Friedman} we have that $w_{\epsilon} \leq u_{\epsilon} \leq v_{\epsilon}$. Moreover from \cite{Amadori}, $\{w_{\epsilon}\}_{\epsilon > 0}$ and $\{v_{\epsilon}\}_{\epsilon > 0}$ are uniformly bounded in 
  $L^{\infty}_{\mathrm{loc}}(\Pi)$.
  Define the following two functions:
  \begin{equation*}
     W := \text{inf}_{\epsilon > 0} \, w_{\epsilon} \; \mathrm{and} \; V := \text{sup}_{\epsilon > 0} \, v_{\epsilon}.
  \end{equation*}
  Clearly therefore $W, V \in L^{\infty}_{\mathrm{loc}}(\Pi)$ and they satisfy the inequality 
  \begin{equation*}
   W \leq w_{\epsilon} \leq u_{\epsilon} \leq v_{\epsilon} \leq V.
  \end{equation*}
This proves \eqref{eq:uLinfinitylocbound}.

Now since $|P_{\epsilon}(t,x)| = \Big|\int_{0}^{x} u_{\epsilon}(t,y) dy\Big| \leq \int_{0}^{x}|u_{\epsilon}(t,y)| dy$, \eqref{eq:PLinfinitylocbound} follows from \eqref{eq:uLinfinitylocbound}. This completes the proof.
\end{proof}
\end{lemma}

\section{Proof of the Main Theorem}\label{sec:MainResult}
In this section we prove Theorem \ref{thm:wellposedness}. Using the compensated compactness method, (see \cite{Tartar, Panov}) we are going to construct a solution of \eqref{eq:OHII} or equivalently of \eqref{eq:OHIII} by passing to the limit in sequence $\{u_{\epsilon}\}_{\epsilon > 0}$ of the viscosity approximations \eqref{eq:OHIV}. The compensated compactness method due to Panov (see Theorem $5$ of \cite{Panov}, or Lemma $2.2$ of \cite{GiuseppeComp}) to be used here can be stated as the following lemma
\begin{lemma}\label{lem:CompensatedCompactness}
 Let $\{v_{\epsilon}\}_{\epsilon > 0}$ be a family of functions defined on $\Pi$. If $\{v_{\epsilon}\}_{\epsilon > 0}$ is uniformly
bounded in $L^{\infty}_{\mathrm{loc}}(\Pi)$ and the family $\{\partial_{t} \eta(v_{\epsilon}) + \partial_{x} 
q(x,v_{\epsilon})\}_{\epsilon > 0}$ is compact in $H^{-1}_{\mathrm{loc}}(\Pi)$ for every convex $\eta \in C^{2}(\R)$, 
where $q_{u}(x,u) = \eta'(u)f_{u}(x,u)$. Then there exist a sequence $\{\epsilon_{k}\}_{k \in \N} \subset \R_{+}$, $\epsilon_{k} \to 0$ as $k \to \infty$, and a map
$v \in L^{\infty}_{\mathrm{loc}}(\Pi)$ such that $v_{\epsilon_{k}} \longrightarrow v$ a.e. and in $L^{p}(\Pi)$ $1 \leq p < \infty$, as $k \to \infty$.
\end{lemma}
The following compact embedding result of Murat \cite{Murat} will be also used,
\begin{lemma}\label{lem:MuratLemma}
 Let $\Omega$ be a bounded open set of $\R^{N}$, $N \geq 2$. Suppose that the sequence $\{\mathcal{L}_{\epsilon}\}_{\epsilon \in \N}$ 
 of distributions is bounded in $W^{-1, \infty}(\Omega)$. In addition, suppose that $\mathcal{L}_{\epsilon} = \mathcal{L}_{1,\epsilon} +
 \mathcal{L}_{2,\epsilon}$; where $\{\mathcal{L}_{1,\epsilon}\}_{\epsilon \in \N}$ lies in a compact subset of $H^{-1}_{\mathrm{loc}}(\Omega)$
 and $\{\mathcal{L}_{2,\epsilon}\}_{\epsilon \in \N}$ lies in a bounded subset of $L^{1}_{\mathrm{loc}}(\Omega)$.
 Then $\{\mathcal{L}_{\epsilon}\}_{\epsilon \in \N}$ lies in a compact subset of $H^{-1}_{\mathrm{loc}}(\Omega)$.
\end{lemma}
First we are going to extract a limit function $u$ from the collection $u_{\epsilon}$ and then we are going to show that this $u$ satisfies \eqref{eq:StrongTraceEntropy}.
\begin{lemma}\label{lem:compencompact}
 The family $\{u_{\epsilon}\}_{\epsilon > 0}$ has a subsequence $\{u_{\epsilon_k}\}_{k \in \N}$ and a limit function
 $u \in L^{\infty}_{\mathrm{loc}}(\Pi)$ such that
 \begin{equation}\label{eq:a.e.convergenceu}
  u_{\epsilon_k} \to u \quad \text{a.e. and in} \quad L^{p}_{\mathrm{loc}}(\Pi), 1 \leq p < \infty.
 \end{equation}
Moreover, we have
  \begin{equation}\label{eq:a.e.convergenceP}
   P_{\epsilon_k} \to P \quad \text{a.e. and in} \quad L^{p}_{\mathrm{loc}}(\R_{+};W^{1,p}_{\mathrm{loc}}(\R_{+})), 1 \leq p < \infty,
  \end{equation}
where 
  \begin{equation*}
    P(t,x) = \int\displaylimits_{0}^{x}u(t,y) dy, \quad t \geq 0, x \geq 0.
  \end{equation*}
  Moreover, \eqref{eq:StrongTraceEntropy} is satisfied.
\begin{proof}
Let . Multiplying the equation \eqref{eq:OHV} by $\eta'(u_{\epsilon})$, we get
\begin{equation*}
 \partial_{t} u_{\epsilon} \eta'(u_{\epsilon}) + \partial_{x} f(x,u_{\epsilon}) \eta'(u_{\epsilon}) = P_{\epsilon} \eta'(u_{\epsilon}) 
 + \epsilon \partial_{xx}^{2}u_{\epsilon}\eta'(u_{\epsilon}),
\end{equation*}
which can be rewritten as
\begin{equation*}
 \partial_{t}\eta(u_{\epsilon}) + f_{u_{\epsilon}}(x,u_{\epsilon}) (u_{\epsilon})_{x}\eta'(u_{\epsilon}) + f_{x}(x,u_{\epsilon})\eta'(u_{\epsilon})
 = P_{\epsilon} \eta'(u_{\epsilon}) + \epsilon \partial_{xx}^{2}u_{\epsilon}\eta'(u_{\epsilon}).
\end{equation*}
From the definition of $q(x,u_{\epsilon})$ we have $q_{u_{\epsilon}}(x,u_{\epsilon}) = \eta'(u_{\epsilon}) f_{u_{\epsilon}}(x, u_{\epsilon})$.
Inserting this into the above expression we get
\begin{equation*}
 \partial_{t}\eta(u_{\epsilon}) + \partial_{x}q(x, u_{\epsilon}) + f_{x}(x,u_{\epsilon})\eta'(u_{\epsilon}) - 
 q_{x}(x, u_{\epsilon}) = P_{\epsilon} \eta'(u_{\epsilon}) + \epsilon 
 \partial_{xx}^{2}u_{\epsilon}\eta'(u_{\epsilon}).
\end{equation*}
This can be written as
\begin{align}\label{eq:EntropyInequalityI}
 \partial_{t}\eta(u_{\epsilon}) + \partial_{x}q(x, u_{\epsilon}) & = \underbrace{\epsilon \partial_{xx}^{2} \eta(u_{\epsilon})}_{\mathcal{L}^{1}_{\epsilon}} - \underbrace{ \epsilon \eta''(u_{\epsilon})(\partial_{x}u_{\epsilon})^{2}}_{\mathcal{L}^{2}_{\epsilon}} + \underbrace{\eta'(u_{\epsilon}) P_{\epsilon}}_{\mathcal{L}^{3}_{\epsilon}}  \nonumber \\ & \quad + \underbrace{ q_{x}(x, u_{\epsilon})}_{\mathcal{L}^{4}_{\epsilon}} - \underbrace{ f_{x}(x,u_{\epsilon})\eta'(u_{\epsilon})}_{\mathcal{L}^{5}_{\epsilon}}.
\end{align}
From Lemma \ref{lm:zeromeanandenergybound} we have
\begin{align*}
 \mathcal{L}^{1}_{\epsilon} \to 0, \; \mathrm{in} \; H^{-1}_{\mathrm{loc}}(\Pi), \; \{\mathcal{L}^{2}_{\epsilon}\}_{\epsilon > 0} \; \mathrm{is} \; \mathrm{uniformly} \; \mathrm{bounded} \; \mathrm{in} \; L^{1}_{\mathrm{loc}}(\Pi).
\end{align*}
To show $\{\mathcal{L}^{3}_{\epsilon}\}_{\epsilon > 0}$ is uniformly bounded in $L^{1}_{\mathrm{loc}}(\Pi)$, let $K$ be any bounded subset of $\Pi$. Then, by Lemma \ref{lm:LocallyBoundLemma},
\begin{equation*}
\|\eta'(u_{\epsilon})P_{\epsilon}\|_{L^{1}(K)} \leq \|\eta'(u_{\epsilon})\|_{L^{\infty}(K)} \|P_{\epsilon}\|_{L^{\infty}(K)} |K|.
\end{equation*}
So it remains to show that $f_{x}(x,u_{\epsilon})\eta'(u_{\epsilon})$ and $q_{x}(x, u_{\epsilon})$ are uniformly bounded in $L^{1}_{\mathrm{loc}}(\Pi)$. To that end observe that
\begin{align*}
 \|f_{x}(x,u_{\epsilon})\eta'(u_{\epsilon})\|_{L^{1}(K)} &= \int\displaylimits_{K}|f_{x}(x,u_{\epsilon})\eta'(u_{\epsilon})| dx dt \, (\mathrm{by} \; \ref{cond:boundedderivative} \; \mathrm{and} \; |\eta'(u_{\epsilon})| \leq C|u_{\epsilon}| )\\
 &\leq \tilde{C} \int\displaylimits_{K} |u_{\epsilon}|^{2} dx dt \, (\mathrm{for} \; \mathrm{some} \; \mathrm{constant} \; \tilde{C}>0)\\
 &< \infty
\end{align*}
So $\{f_{x}(x,u_{\epsilon})\eta'(u_{\epsilon})\}_{\epsilon > 0}$ is uniformly bounded in $L^{1}_{\mathrm{loc}}(\Pi)$.
Similarly we have 
\begin{align*}
 \|q_{x}(x,u_{\epsilon})\|_{L^{1}(K)} &= \int\displaylimits_{K} \Big|\int\displaylimits_{0}^{u_{\epsilon}} \eta'(v)f_{xv}(v) dv\Big|dxdt \, (\mathrm{by} \; \ref{cond:boundedderivative}) \\
 &\leq C \int\displaylimits_{K}\int_{0}^{u_{\epsilon}}|\eta'(v)| dv dx dt \\ 
 &< \infty.
\end{align*}
Consequently, $\{q_{x}(x,u_{\epsilon})\}_{\epsilon > 0}$ is uniformly bounded in $L^{1}_{\mathrm{loc}}(\Pi)$.

Therefore, by Lemma \ref{lem:MuratLemma} we can conlcude that
\begin{equation}\label{eq:compactbound}
 \{\partial_{t} \eta(u_{\epsilon}) + \partial_{x} q(x, u_{\epsilon})\}_{\epsilon > 0} \quad \text{lies in a compact subset of} \quad 
 H^{-1}_{\mathrm{loc}}(\Pi).
\end{equation}
Therefore using the $L^{\infty}_{\mathrm{loc}}$ bound obtained from Lemma \ref{lm:LocallyBoundLemma}, \eqref{eq:compactbound} and 
Lemma \ref{lem:CompensatedCompactness} we can conclude that there exists a subsequence $\{u_{\epsilon_{k}}\}_{k \in \N}$ and
a limit function $u \in L^{\infty}_{\mathrm{loc}}(\Pi)$ such that \eqref{eq:a.e.convergenceu} holds.
By the H\"{o}lder inequality and the definition of $P_{\epsilon}$, \eqref{eq:a.e.convergenceP} follows from \eqref{eq:a.e.convergenceu}.

We remark that the entropy inequality \eqref{eq:entropyinequality} can be obtained from \eqref{eq:EntropyInequalityI} by the standard argument of letting $\epsilon \to 0$ and using convexity of $\eta(\cdot)$. Thus by \cite[Theorem $1.1$]{GiuseppeY}, strong trace $u^{\tau}_{0}$ for $u$ on $x=0$ does exist. Now we are going to prove \eqref{eq:StrongTraceEntropy}. From the Definition \ref{def:EntropySoln} for \eqref{eq:OHV} and using \eqref{eq:entropyinequality} we get for Kr\v{u}zkov entropy/entropy flux pair $(\eta,q)$
\begin{align*}
& \partial_{t}|u_{\epsilon_{k}} - c| + \partial_{x}\Big(\mathrm{sign}(u-c)(f(x,u) - f(x,c))\Big) \\ & \qquad - \mathrm{sign}(u_{\epsilon_{k}} - c) P_{\epsilon_{k}} - \epsilon_{k} \partial_{xx}^{2} |u_{\epsilon_{k}} - c| \leq 0.
\end{align*}
Multiplying by a non-negative test function $\phi \in C^{2}_{c}(\Pi)$ and integrating over $\Pi$, we get
\begin{align*}
& \iint\displaylimits_{\Pi} \Big(|u_{\epsilon_{k}} - c|\partial_{t} \phi + \Big(\mathrm{sign}(u_{\epsilon_{k}}-c)(f(x,u) - f(x,c))\Big) \partial_{x} \phi \nonumber \\ & \qquad \qquad - \mathrm{sign}(u_{\epsilon_{k}}-c) f_{x}(x,c) \phi + \mathrm{sign}(u_{\epsilon_{k}} -c) P_{\epsilon_{k}} \phi \Big) dt \, dx \nonumber \\ & \qquad \qquad - \epsilon_{k} \iint\displaylimits_{\Pi} \partial_{x} \, |u_{\epsilon_{k}} - c| \, \partial_{x} \, \phi \, dt \, dx + \int\displaylimits_{\R_{+}} |u_{0}(x) - c| \phi(x,0) dx \nonumber \\
& \qquad \qquad + \int\displaylimits_{\R_{+}} \mathrm{sign}(c) f(0,c) \phi(0,t) dt - \epsilon_{k} \int\displaylimits_{\R_{+}} \partial_{x}|u_{\epsilon_{k}}(0,t) - c| \phi(0,t) dt \geq 0.
\end{align*}
Invoking Lemmas \ref{lm:zeromeanandenergybound}, \ref{lm:LocallyBoundLemma}, and \ref{lem:compencompact}, letting $k \to \infty$, we have
\begin{align}\label{eq:StrongTraceEntropyII}
& \iint\displaylimits_{\Pi} \Big(|u - c|\partial_{t} \phi + \Big(\mathrm{sign}(u-c)(f(x,u) - f(x,c))\Big) \partial_{x} \phi \nonumber \\ & \qquad \qquad - \mathrm{sign}(u-c) f_{x}(x,c) \phi + \mathrm{sign}(u-c) P \phi \Big) dt \, dx \nonumber \\ & \qquad \qquad + \int\displaylimits_{\R_{+}} |u_{0}(x) - c| \phi(x,0) dx  + \int\displaylimits_{\R_{+}} \mathrm{sign}(c) f(0,c) \phi(0,t) dt \\ & \qquad \qquad - \mathrm{lim}_{k \to \infty} \epsilon_{k} \int\displaylimits_{\R_{+}} \partial_{x}|u_{\epsilon_{k}}(0,t) - c| \phi(0,t) dt\geq 0. \nonumber
\end{align}
Consequently to show \eqref{eq:StrongTraceEntropy} it is enough to prove that
\begin{align}\label{eq:StrongTraceEntropyI}
\begin{split}
& \mathrm{lim}_{k \to \infty} \epsilon_{k} \int\displaylimits_{\R_{+}} \partial_{x}|u_{\epsilon_{k}}(0,t) - c| \phi(0,t) dt = \int\displaylimits_{\R_{+}} \mathrm{sign}(c) f(0,u^{\tau}_{0}(t)) \phi(0,t) dt.
\end{split}
\end{align}
In order to prove this we need to employ a particular choice of test function. Let $\{\Psi_{m}\}_{m \in \N} \subset C^{\infty}_{c}(\R)$ be a sequence of non-negative test functions satisfying
\begin{equation}
\begin{cases}
\Psi_{m}(0) = 1, \; \mathrm{for} \; \mathrm{all} \; m \in \N, \\
|\Psi'_{m}| \leq m, \; \mathrm{and}\\
\Psi_{m}(x) = 0, \; \mathrm{for} \; \mathrm{all} \; x \geq \frac{1}{m} .
\end{cases}
\end{equation}
Multiplying the equation \eqref{eq:OHV} by the test function $\Psi_{m}(x) \phi(x,t)$ we get after an integration by parts
\begin{align}\label{eq:EntropyTrace}
& \iint\displaylimits_{\Pi} \Big(u_{\epsilon_{k}} \partial_{t} \phi \Psi_{m} + f(x,u_{\epsilon_{k}}) (\Psi_{m}\partial_{x}\phi + \Psi'_{m} \phi) + P_{\epsilon_{k}} \Psi_{m} \phi\Big) \, dt \, dx \nonumber \\ & \qquad - \iint\displaylimits_{\Pi} \epsilon_{k} \partial_{x} u_{\epsilon_{k}} (\Psi_{m}\partial_{x}\phi + \Psi'_{m} \phi) \, dt \, dx + \int\displaylimits_{\R_{+}} u_{0}(x) \phi(x,0) \Psi_{m}(x) \, dx \\ & \qquad - \int\displaylimits_{\R_{+}} f(0,u_{\epsilon_{k}(0,t)}) \phi(0,t) \, dt  - \epsilon_{k} \int\displaylimits_{\R_{+}} \partial_{x} u_{\epsilon_{k}}(0,t) \phi(0,t) \, dt= 0 \nonumber.
\end{align}
Employing the strong convergence $u_{\epsilon_{k}} \to u$ from Lemma \ref{lem:compencompact}, passing to the limit $k \to \infty$, $m \to \infty$ respectively and using the properties of $\Psi_{m}$ in the above relation \eqref{eq:EntropyTrace} we get
\begin{equation*}
\mathrm{lim}_{k \to \infty} \epsilon_{k} \int\displaylimits_{\R_{+}} \partial_{x} u_{\epsilon_{k}}(0,t) \phi(0,t) dt = - \int\displaylimits_{\R_{+}} f(0,u^{\tau}_{0}(t)) \phi(0,t) dt,
\end{equation*}
which in turn proves \eqref{eq:StrongTraceEntropyI}. Combining \eqref{eq:StrongTraceEntropyII} and \eqref{eq:StrongTraceEntropyI} we have obtained the desired inequality \eqref{eq:StrongTraceEntropy}.

This completes the proof.
\end{proof}
\end{lemma}  
Consequently we have established the existence of an entropy solution (in the sense of Definition \ref{def:EntropySoln}) $u(x,t)$ of the equation \eqref{eq:OHII} or equivalently of \eqref{eq:OHIII}. Now in order to prove the uniqueness of entropy solutions we are going to prove \eqref{eq:contractivity}, \textit{i.e.} we will prove Theorem \ref{thm:wellposedness}.

\begin{proof}{(of Theorem \ref{thm:wellposedness} :)}
Let $u$ and $v$ be two entropy solutions of \eqref{eq:OHII} or equivalently \eqref{eq:OHIII}. We will use the doubling of variables. For $\Pi := (0,\infty)^{2}$ and $\Pi^{2} := (0,\infty)^{4}$ let $\phi(t,\tau,x,y) \in C^{\infty}_{c}(\Pi^{2})$ be a non-negative test function. Since $u$ and $v$ are entropy solutions of
\eqref{eq:OHIII}, we have 
\begin{align}\label{eq:doubling1}
 & \iint\displaylimits_{\Pi} \Big[|u(x,t) - v(y,\tau)|\partial_{t}\phi(t,\tau,x,y) + [f(x,u(x,t)) - f(y,v(y,\tau))] \nonumber \\ & \qquad \qquad  \mathrm{sign}(u(x,t) - v(y,\tau))\partial_{x}\phi(t,\tau,x,y) \nonumber \\ & \qquad \qquad -  \mathrm{sign}(u(x,t) - v(y,\tau))[f_{x}(x,v(y,\tau)) - P_{u}(x,t)]\phi(t,\tau,x,y) \Big] dt dx \geq 0, 
\end{align}
and
\begin{align}\label{eq:doubling2}
 & \iint\displaylimits_{\Pi} \Big[|v(y,\tau) - u(x,t)|\partial_{\tau}\phi(t,\tau,x,y) + [f(y,v(y,\tau)) - f(x,u(x,t))] \nonumber \\ & \qquad \qquad \mathrm{sign}(v(y,\tau) - u(x,t))\partial_{y}\phi(t,\tau,x,y) \nonumber \\ & \qquad \qquad - \mathrm{sign}(v(y,\tau) - u(x,t))[f_{y}(y,u(x,t)) - P_{v}(y,\tau)]\phi(t,\tau,x,y) \Big] d\tau dy \geq 0.
\end{align}
Then integrating \eqref{eq:doubling1} with respect to $\tau$, $y$; \eqref{eq:doubling2} with respect to $t$, $x$; and 
adding the two outcomes we obtain,
\begin{align}\label{eq:doubling3}
 &\iiiint\displaylimits_{\Pi^{2}}  \Big[|u(x,t) - v(y,\tau)|(\partial_{t}\phi(t,\tau,x,y) + \partial_{\tau}\phi(t,\tau,x,y)) + [f(x,u(x,t)) \nonumber \\ & \qquad \qquad - f(y,v(y,\tau))] \mathrm{sign}(u(x,t) - v(y,\tau))(\partial_{x}\phi(t,\tau,x,y) + \partial_{y}\phi(t,\tau,x,y)) \nonumber \\ & \qquad \qquad + \mathrm{sign} (u(x,t) - v(y,\tau))(P_{u}(x,t) - P_{v}(y,\tau))\phi(t,\tau,x,y) \nonumber \\ & \qquad \qquad - \mathrm{sign}(u(x,t) - v(y,\tau)) (f_{x}(x,v(y,\tau)) - f_{y}(y,u(x,t))) \phi(t,\tau,x,y)\Big] dt d\tau dx dy \geq 0.
\end{align}
For $\rho_{\epsilon} \to \delta_{0}$ as $\epsilon \to 0$, where $\delta_{0}$ is the Dirac mass concentrated at $0$, where 
\begin{equation}\label{eq:mollifier}
 \rho_{\epsilon}(z) := \epsilon \rho(\epsilon z), \quad \text{and} \quad \alpha_{\epsilon}(z) := \int_{-\infty}^{z}\rho_{\epsilon}(x)dx,
\end{equation}
for some non-negative $\rho \in C^{\infty}_{c}([-1,1])$ with total mass being $1$. Now let us define the particular test function
\begin{equation}\label{eq:particulartestfuncn}
 \phi_{\epsilon}(t,\tau,x,y) = \psi\Big(\frac{t + \tau}{2},\frac{x + y}{2}\Big) \rho_{\epsilon}\Big(\frac{\tau - t}{2}\Big) 
 \rho_{\epsilon}\Big(\frac{y - x}{2}\Big),
\end{equation}
where $\psi \in C^{\infty}_{c}(\Pi)$ is a non-negative, test function. Inserting the function \eqref{eq:particulartestfuncn} into the last inequality \eqref{eq:doubling3}, we get
\begin{align}\label{eq:doubling4}
 &\iiiint\displaylimits_{\Pi^{2}} \Big[\rho_{\epsilon}\Big(\frac{\tau - t}{2}\Big) \rho_{\epsilon}
 \Big(\frac{y - x}{2}\Big)\Big\{|u(t,x) - v(\tau,y)|\partial_{t}\psi\Big(\frac{t + \tau}{2},\frac{x + y}{2}\Big) \nonumber \\
 &\qquad \qquad + \Big(f(x,u(t,x)) - f(y,v(\tau,y))\Big)\mathrm{sign}(u(t,x) - v(\tau,y))\partial_{x}\psi\Big(\frac{t + \tau}{2},\frac{x + y}{2}\Big)
 \Big\} \nonumber \\ &\qquad \qquad + \gamma \psi\Big(\frac{t + \tau}{2},\frac{x + y}{2}\Big)\rho_{\epsilon}\Big(\frac{\tau - t}{2}\Big) \rho_{\epsilon}
 \Big(\frac{y - x}{2}\Big) \mathrm{sign} (u(t,x) - v(\tau,y)) \nonumber \\ &\qquad \qquad \qquad (P_{u}(t,x) - P_{v}(\tau,y))
 - \mathrm{sign}(u(t,x) - v(\tau,y)) (f_{x}(x,v(\tau,y)) \nonumber \\ &\qquad \qquad - f_{y}(y,u(t,x))) 
 \psi\Big(\frac{t + \tau}{2},\frac{x + y}{2}\Big)\rho_{\epsilon}\Big(\frac{\tau - t}{2}\Big) \rho_{\epsilon}
 \Big(\frac{y - x}{2}\Big)\Big]dt d\tau dx dy \geq 0.
\end{align}
By standard limiting argument of doubling of variable technique, passing to the limit as $\epsilon \to 0$ we obtain from the previous inequality \eqref{eq:doubling4} that for all test functions $\psi$ as mentioned above,
\begin{align}\label{eq:doubling5}
 &\iint\displaylimits_{\Pi} \Big[|u(t,x) - v(t,x)| \partial_{t} \psi + \Big(f(x,u) - f(x,v)) \mathrm{sign}(u(t,x) - v(t,x))\Big)
 \partial_{x}\psi \Big] dt dx \nonumber \\ & \qquad \qquad + \iint\displaylimits_{\Pi} \mathrm{sign}(u(t,x) - v(t,x)) ((P_{u}(t,x) - P_{v}(t,x))
 \psi dt dx \nonumber \\ &\qquad \qquad + \iint\displaylimits_{\Pi} \mathrm{sign}(u(t,x) - v(t,x)) \Big(f_{x}(x,u(t,x)) - f_{x}(x,v(t,x))\Big) \psi dt dx \geq 0.
\end{align}
Following Kru\v{z}kov's argument \cite{Kruzkov} if we consider the sets for $T$, $R > 0$
\begin{equation}
 \Omega_{R,T} := \{(t,x) \in [0,T] \times [0,R] ; \quad 0\leq s \leq t , \quad 0 \leq x \leq R + L(t-s)\}, 
\end{equation}
and define the following non-negative test function
\begin{equation*}
 \phi_{\epsilon}(t,x) := [\alpha_{\epsilon}(s) - \alpha_{\epsilon}(s - t)][1 - \alpha_{\epsilon}(x - R - L(t-s))],
\end{equation*}
where $\alpha_{\epsilon}$ is defined in \eqref{eq:mollifier} and $L$ is defined in \ref{cond:uLipschitz}. Clearly observe that $\phi_{\epsilon}$ is an approximation of the characteristic function of $\Omega_{R,T}$. From definition $\alpha_{\epsilon}' = \rho_{\epsilon} \geq 0$. Using $\phi_{\epsilon}$ as the test function in \eqref{eq:doubling5} and similarly as before letting $\epsilon \to 0$, we get
\begin{align}\label{eq:uniquenessestimate}
 \|u(\cdot,t) - v(\cdot,t)\|_{L^{1}(0,R)} &\leq \|u_{0} - v_{0}\|_{L^{1}(0,R+Lt)} \nonumber \\
 & \qquad + \int\displaylimits_{\Omega_{R,T}} \mathrm{sign} (u(x,t) - v(x,t)) (P_{u} - P_{v}) dx ds \nonumber \\
 & \qquad + \int\displaylimits_{\Omega_{R,T}} \mathrm{sign}(u(x,t) - v(x,t)) \Big(f_{x}(x,u) - f_{x}(x,v)\Big) dx ds .
\end{align}
With
\begin{equation}
I(s) := [0, R+L(t-s)],
\end{equation}
note that
\begin{align}\label{eq:spacedependentextra}
 \int\displaylimits_{\Omega_{R,T}}  \mathrm{sign}(u - v) \Big(f_{x}(x,u) - f_{x}(x,v)\Big) dx ds &\leq \int\displaylimits_{0}^{t} \int\displaylimits_{I(s)} |f_{x}(x,u) - f_{x}(x,v)| dx ds \nonumber \\ &\qquad \qquad (\mathrm{Using} \; \ref{cond:xLipschitz})
 \nonumber \\ &\leq \int\displaylimits_{0}^{t} \int\displaylimits_{I(s)} L_{1} |u - v| dx ds \nonumber \\ &\leq \int\displaylimits_{0}^{t} L_{1}\|u - v\|_{L^{1}(I(s))} ds.
\end{align}
Since
\begin{align}\label{eq:nonlocalestimate}
  \int\displaylimits_{\Omega_{R,T}} \text{sign}(u - v)(P_{u} - P_{v}) ds dx &\leq  \int\displaylimits_{0}^{t} \int\displaylimits_{I(s)} |P_{u} - P_{v}| ds dx \nonumber \\ &\leq  \int\displaylimits_{0}^{t} \int\displaylimits_{I(s)} \Big(\Big|\int\displaylimits_{0}^{x}|u-v| dy\Big|\Big) ds dx \nonumber \\ &\leq  \int\displaylimits_{0}^{t} \int\displaylimits_{I(s)} \Big(\Big|\int_{I(s)}|u-v| dy\Big|\Big) ds dx \nonumber \\
 &= \int\displaylimits_{0}^{t} |I(s)| \; \|u(\cdot,s) - v(\cdot,s)\|_{L^{1}(I(s))} ds,
\end{align}
and,
\begin{equation}\label{eq:I(s)bound}
 |I(s)| = R + L(t-s) \leq R + Lt \leq R + LT.
\end{equation}
We consider the following continuous function:
\begin{equation}
 G(t) := \|u(\cdot,t) - v(\cdot,t)\|_{L^{1}(I(t))}, \quad t \geq 0.
\end{equation}

Then we can combine \eqref{eq:uniquenessestimate}, \eqref{eq:spacedependentextra}, \eqref{eq:nonlocalestimate} and \eqref{eq:I(s)bound} to obtain
\begin{equation}
 G(t) \leq G(0) +  \int_{0}^{t}(|I(s)| + L_{1})G(s) ds \quad \text{with} \quad |I(s)| = R + L(t-s).
\end{equation}

Consequently, by Gronwall's inequality we can conclude:
\begin{equation*}
 G(t) \leq G(0) e^{ \int_{0}^{t} ( |I(s)| + L_{1})ds}, \quad \text{for a.e.} \quad 0 < t < T,
\end{equation*}
i.e.
\begin{equation*}
 G(t) \leq G(0) e^{(Rt+\frac{1}{2}Lt^{2}) + L_{1}t} \leq G(0) e^{(R+\frac{1}{2}LT)t + L_{1}T}, \quad \text{for a.e.} \quad 0 < t < T.
\end{equation*}
Consequently we have the estimate \eqref{eq:contractivity}, namely
\begin{equation}
 \|u(\cdot,t) - v(\cdot,t)\|_{L^{1}((0,R))} \leq e^{Ct}\|u(\cdot,0) - v(\cdot,0)\|_{L^{1}(0,R+Lt)}, \quad \text{for a.e.} \quad 0 < t < T,
\end{equation}
where the constant $C$ depends on $T$, $R$, $L_{1}$ and $L$.

This completes the proof.
\end{proof}

\bibliography{reference} \bibliographystyle{amsplain}

\end{document}